# On the solution of the Graph Isomorphism Problem

# Part II

Leonid I. Malinin (1922-1999)

Natalia L. Malinina

(Moscow Aviation Institute (State University of Aerospace Technologies))

September 15, 2010

**Abstract**

The presented material continues the previous article and also is devoted to the equivalent conversion between the graphs. The examining of the transformation of the vertex graphs into the edge graphs (together with the opposite transformation) illustrates the reasons of the appearance of NP-completeness from the point of view of the graph theory. We suggest that it also illustrates the synchronous possibility and impossibility of the struggle with NP-completeness.

## A converting of the directed graphs

## 1 Introduction

The previous paper was devoted to the operation of the normalization of the arbitrary adjacency matrix and to the operation of the reduction of the received matrix. Very briefly was described the operation of the constructing the $F$ matrix according to the $R$ matrix. But at transforming the vertex graphs into the edge graphs such an operation as the construction of the $F$ matrix with the help of either canonical or quasicanonical $R$ matrix is used. This operation is undoubtedly interesting when constructing the net models of such systems, which complexity may be bind with the cyclomatic number.

A transfer from the $L$ matrix to the $R$ matrix may be presented by the scheme in fig. 1.

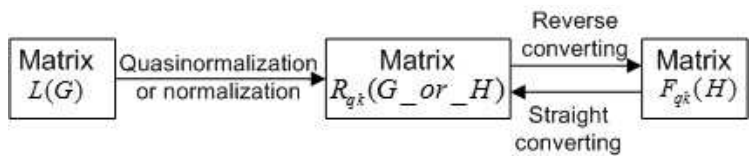

Fig. 1

Where:

$L(G)$ – An adjacency matrix of the initial $G$ graph's vertexes;

$R_{q(k)}(G \text{ or } H)$ – A quasicanonical or canonical dual adjacency matrix (of either the $G$ graph's vertexes or the $H$ graph's edges);

$F_{q(k)}(H)$ – A corresponding to $R_{q(k)}$ (either $G$ or $H$) matrix – either the quasicanonical or the canonical adjacency matrix of the $H$ graph's vertexes;

A transfer from the $R$ matrix to the $F$ matrix represents the transfer from the $G$ graph to the $H$ graph in the matrix's form as the operation of constructing the edge graph by the given vertex graph. Let us denote this operation as the transformation or the converting of the vertex graph into the edge graph.

An operation of constructing the vertex graph by the given edge graph is the opposite operation. It is simple and for the undirected graphs is described in [1]. Let us examine both operations for the directed graphs.

Let us agree to denote the operation of constructing the vertex graph by the given edge graph as the straight converting, and the operation of constructing the edge graph by the given vertex graph – the reverse converting (fig.1).

Both of the operations have the following evident features:

1. Every directed $H$ graph may be subjected to the consecutive straight converting any number of times.

2. The directed $G$ graph may be subjected to the consecutive reverse converting if and only if the adjacency matrix of its vertexes, being considered as the $R$ matrix, has either the quasicanonical or the canonical form.

3. If the $H$ graph is a result of the single straight converting of some $P$ graph, then it may be subjected to the reverse converting at least once. If the $H_j$ graph is a result of the single straight converting of the $H_1$ graph for $(j-1)$ times, then it may be subjected to the operation of the reverse converting at least the $j$ times.

Let us denote the operation of the straight converting as $J$, and the operation of the reverse converting as $D$, that is:

$J: F \to R$

$D: R \to F$

It is evident that $J = D^{-1}$

Let us agree to denote both $F$ and $R$ matrixes and corresponding to them both $H$ and $G$ graphs without either $q$ or $k$ indexes, if there is no need to distinguish the cases of the converting either quasicanonical or canonical graphs.

The $F$ matrix differs from the $R$ matrix. The $r_{ij} = 1$ elements in the $R$ matrix define only the fact of the binary $(q_i < q_j)$ relation of either two $G$ graph's both $q_i$ and $q_j$ vertexes or two $H$ graph's both $q_i$ and $q_j$ edges. On the contrary, in the $F$ matrix the $f_{hg} = 1$ elements on the one hand define the fact of the binary ($v_h < v_g$) relation of the two $H$ graph's both $v_h$ and $v_g$ vertexes, and on the other hand, they have some physical or mathematical content as the operators, which transform the input data of the process's model (vertex $v_h$) into the output data (vertex $v_g$). Thus, the $F$ matrix is denoted as the operator adjacency matrix of the $H$ graph's vertexes.

The $F$ matrix must meet the definite requirements; it follows from its definition and from the approach to its receiving. If the $F$ matrix is received from the quasicanonical $R_q$ matrix, then it is called the quasicanonical, is denoted as $F_q$ and must meet the following requirements:

- All the matrix's elements, situated on the main diagonal, must be equal to zero;
- There could not be more than one empty column (corresponding to the initial vertex) in the $F$ matrix;
- If there is an empty column in the $F$ matrix, then in the corresponding line there must be exactly one non-zero element;
- There could not be more than one empty line (corresponding to the final vertex) in the $F$ matrix;
- If there is an empty line in the $F$ matrix, then in the corresponding column there must be exactly one non-zero element.

If the $F$ matrix is received from the canonical $R_k$ matrix, then it is called the canonical, denoted as $F_k$ and must, in addition to the listed ones, meet the following requirements:

- If there are more than one non-zero element in the $F$ matrix's the $h$ line, then there must be exactly one non-zero element in the $g = h$ column;
- If there are more than one non-zero element in the $F$ matrix's the $g$ column, then there must be exactly one non-zero element in the $h = g$ line.

Any square matrix, which meets the listed above requirements, may be adopted accordingly as either $F_q$ or $F_k$ matrix and can be transformed into either $R_q$ or $R_k$ matrix.

Since the $F$ matrix can be either the quasicanonical one or the canonical one, then we'll distinguish either the quasicanonical or the canonical converting. It follows from theorems 1 and 4 [2] that at the canonical straight converting the cyclomatic number remains constant, and at the quasicanonical one – it increases.

Let us point to one peculiarity of the straight converting. At the straight converting every edge of some $H_1$ graph, accepted as the edge graph, is transformed into the vertex of new $H_2$ graph. If the $H_1$ graph had both initial and final edges, then the $H_2$ graph will have both initial and final vertexes. If we decide subjecting the $H_2$ graph to the operation of the straight converting once more, then it must be completed with new both initial and final edges. In this case it is considered as the edge graph, and according to its definition, must have the edges as both initial and final elements. So, if the initial $H_1$ graph has either initial or final edges (or both concurrently), then the graphs, obtained as a result of the consecutive application of the straight converting operation, should on every step of the converting be increased with either initial or final edges or by both of them concurrently.

At the indicated condition the straight converting operation can be done in series any $j$ times. So, let us write the following expressions:

$$\left. \begin{array}{l} R_1 = J_1(F_1) \\ R_2 = J_1(J_1(F_1)) = J_2(F_1) = J_1(R_1) \\ \ldots \ldots \\ R_j = J_1(R_{j-1}) = J_j(F_1) \end{array} \right\}$$

Where:

$j = 1,2 \ldots M$, $M$ – an arbitrary large number;

$(j-1)$ – A number of the operation of the straight converting, which is done in series;

$R_j$ - A $R$ matrix, obtained as a result of the $F_1$ matrix's straight converting being done the $(j-1)$ number of times.

For the operation of the reverse converting:

$$\left. \begin{array}{l} F_1 = D_1(R_1) \\ F_2 = D_1(D_1(R_1)) = D_2(R_1) = D_1(F_1) \\ \ldots \ldots \\ F_\alpha = D_1(F_{(\alpha-1)}) = D_\alpha(R_1) \end{array} \right\}$$

Where:

$\alpha = 1,2, \ldots \mu$

$(\alpha - 1)$ – A number of the reverse converting operation done in series;

$F_\alpha$ – A $F$ matrix, which is obtained as a result of the $R_1$ matrix's reverse converting being done $(\alpha - 1)$ times.

A $\mu$ value could not be the arbitrary one. A top boundary of the $\mu$ value is defined by the initial for reverse converting $R_1$ matrix. A review of the reverse converting is possible till newly received both $F_\alpha$ and the $R$ matrixes still meet the requirements of theorems 1 and 4 [2].

## 2   An example of the straight converting

Let us take up the example of the straight converting of the arbitrary directed graph. The edge $H_{q1}$ graph, presented in fig.2, is both given and accepted for the initial graph.

A corresponding $F_{q1}$ matrix is presented in fig. 3. Table 1 (fig. 3) is arranged according to the $F_{q1}$ matrix. An arrangement is evident from the $F_{q1}$ matrix and table 1 (fig. 3), and does not need any explanations.

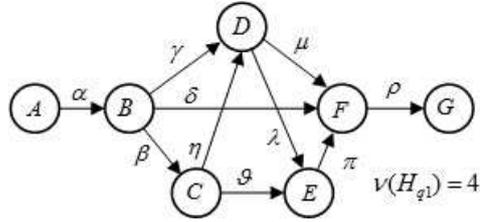

Fig. 2.

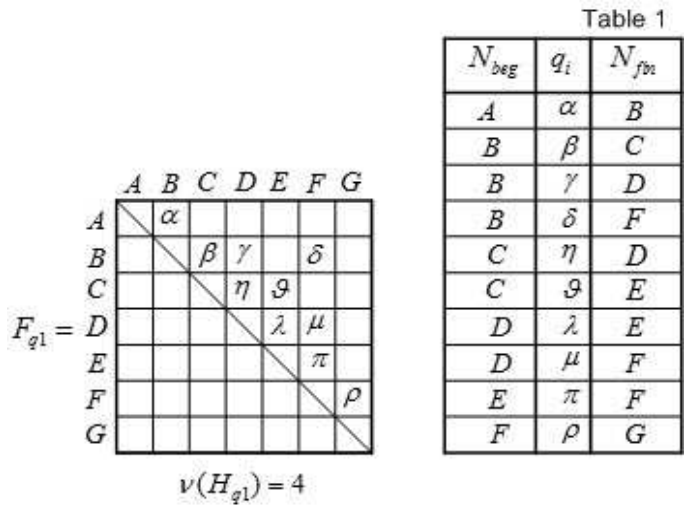

Fig. 3.

Table 2 (fig 4) is arranged accordingly to table 1(fig. 3). A difference between table 1 and table 2 is in the fact that both the first and the second columns of table 1 are joined into the second column of table 2.

Thus in table 2 an ordered pair of the matrix's rows corresponds to every $f_{hg} = i$ element of the $F$ matrix. The $R_{q1}$ matrix is arranged the way you can see below (fig. 4).

$$r_{ij} = \begin{cases} 1, \text{ if the second element of the ordered pair } (N_{beg}, N_{fin}) \\ \quad \text{from the } i - \text{line of the table 2 is congruent to the first} \\ \quad \text{element of the ordered pair } (N_{beg}, N_{fin}) \text{ from the} \\ \quad j - \text{line of the table 2} \\ 0 - \text{otherwise} \end{cases}$$

Later the $R_{q1}$ matrix is transformed into the operational $F_{k2}$ matrix: the $r_{ij} = 1$ elements in the $R_{q1}$ matrix are replaced by the arbitrary symbols (for example, $a, b, c$ ...and so on) and are added, if it is necessary, one both initial and final row (the initial and the final elements of the basic set of the new $H_{k2}$ graph).

A $F_{k2}$ matrix for our example is presented in fig.5.

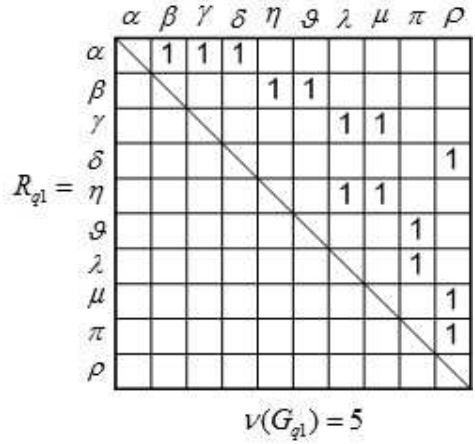

| Vertexes of graph $G_{q1}$ | Pairs of vertexes of graph $H_{q1}$ |
|---|---|
| α | AB |
| β | BC |
| γ | BD |
| δ | BF |
| η | CD |
| ϑ | CE |
| λ | DE |
| μ | DF |
| π | EF |
| ρ | FG |

Table 2

Fig. 4.

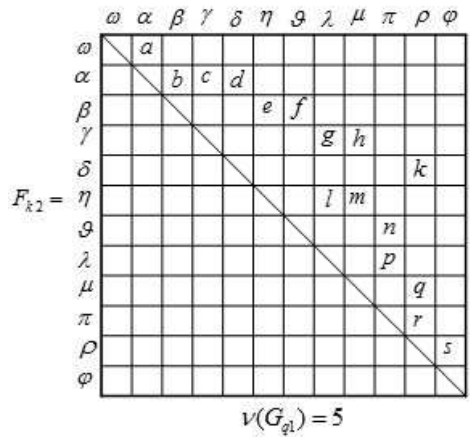

| $N_{beg}$ | $q_i$ | $N_{fin}$ |
|---|---|---|
| ω | a | α |
| α | b | β |
| α | c | γ |
| α | d | δ |
| β | e | η |
| β | f | ϑ |
| γ | g | λ |
| γ | h | μ |
| δ | k | ρ |
| η | l | λ |
| η | m | μ |
| ϑ | n | π |
| λ | p | π |
| μ | q | ρ |
| π | r | ρ |
| ρ | s | φ |

Table 3

Fig. 5.

A matrix is enlarged with the initial ω row and the final φ row. Accordingly both $f_{\omega\alpha} = a$ and $f_{\rho\varphi} = s$ elements are added to the matrix. A $H_{k2}$ graph is constructed by the $F_{k2}$ matrix (fig. 6).

A graph appears to be the canonical one, because it has none complicated vertexes. Its cyclomatic $\nu(H_{k2})$ number is equal to 5, while the cyclomatic $\nu(H_{q1})$ number of the $H_{q1}$ graph is equal to 4.

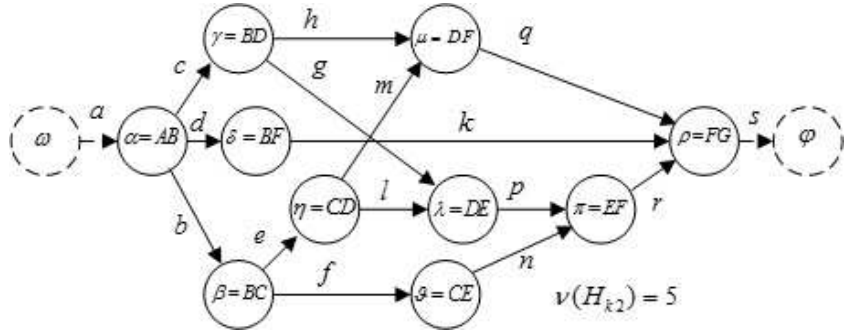

Fig. 6.

Let us note that every $H_{k2}$ graph's vertex corresponds to some vertex's ordered pair of the $H_{q1}$ graph. It is evident that there are as many vertexes in the $H_{k2}$ graph as there are the ordered vertex's pairs in the $H_{q1}$ graph (we do not take into account both added $\omega$ and $\varphi$ vertexes). Equity of this statement follows from the approach of constructing the $H_{k2}$ graph by the $H_{q1}$ graph.

Let's execute the converting operation once more. For this we'll arrange table 3 (fig. 5) by the $F_{k2}$ matrix similarly like we arranged table 1.

The first two columns of table 4 are arranged by table 3 (similarly to table 2). Using table 2, let's decode the ordered pairs of the $F_{k2}$ matrix's rows, which are written in the second column of table 4 (fig. 7).

Table 4

| | $G_{q2}$ graph vertexes | Pairs of $H_{q2}$ graph vertexes | Pairs of towes of $H_{q1}$ graph vertexes | | Threes of $H_{q1}$ graph vertexes | | |
|---|---|---|---|---|---|---|---|
| a | $\omega\alpha$ | $\omega$ | AB | $\omega$ | A | B | |
| b | $\alpha\beta$ | AB | BC | A | B | C | |
| c | $\alpha\gamma$ | AB | BD | A | B | D | |
| d | $\alpha\delta$ | AB | BF | A | B | F | |
| e | $\beta\eta$ | BC | CD | B | C | D | |
| f | $\beta\vartheta$ | BC | CE | B | C | E | |
| g | $\gamma\lambda$ | BD | DE | B | D | E | |
| h | $\gamma\mu$ | BD | DF | B | D | F | |
| k | $\delta\rho$ | BF | FG | B | F | G | |
| l | $\eta\lambda$ | CD | DE | C | D | E | |
| m | $\eta\mu$ | CD | DF | C | D | F | |
| n | $\vartheta\pi$ | CE | EF | C | E | F | |
| p | $\lambda\pi$ | DE | EF | D | E | F | |
| q | $\mu\rho$ | DF | FG | D | F | G | |
| r | $\pi\rho$ | EF | FG | E | F | G | |
| s | $\rho\varphi$ | FG | $\varphi$ | F | G | $\varphi$ | |

Fig. 7.

For example, the $\alpha\beta$ pair (the second line, the second column of table 4) includes both $\alpha = AB$ and $\beta = BC$ elements. Therefore, in the second line of the third column of table 4 the pairs of deuces of both $AB$ and $BC$ elements are written. All the third column of table 4 is filled up in such a way.

Let's notice that in all lines of the third column of table 4 the second elements of the first deuces are congruent to the first elements of the second deuces.

Therefore, every ordered pair of the deuces corresponds to the ordered triple of the $H_{q1}$ graph's vertexes, excluding both the first and the final ordered triples, each including two vertexes of the initial graph and one vertex, added to the $H_{k2}$ graph. For example, the pair of both deuces: $AB$ and $BC$ corresponds to the $ABC$ triple.

The ordered triples are written in the fourth column of table 4. The $R_{k2}$ matrix is arranged (fig. 7) according to this column of table 4. Then the $R_{k2}$ matrix is transformed into the $F_{k3}$ matrix (fig. 8) similarly to the $F_{k2}$ matrix was arranged.

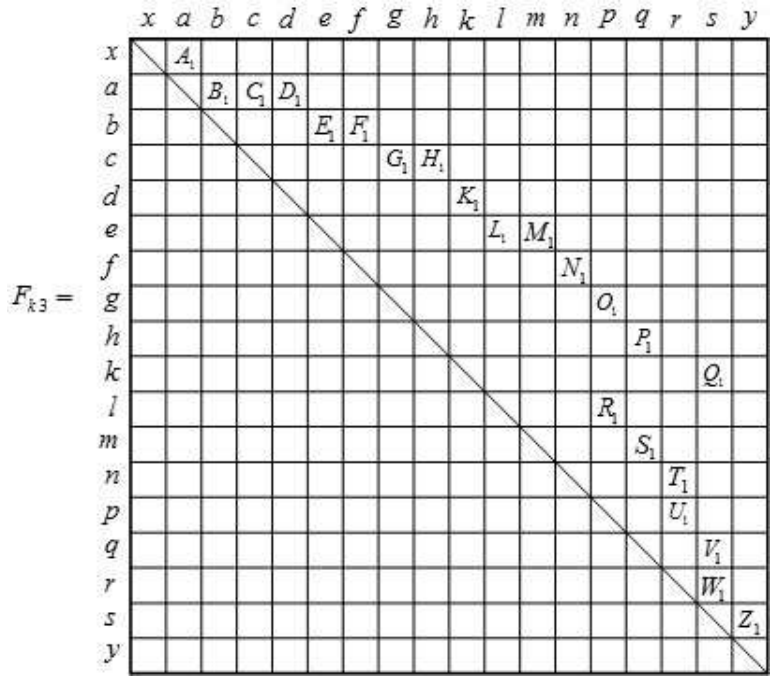

Fig. 8.

Then the $H_{k3}$ graph is constructed by this matrix. It appears to be the canonical graph. Its cyclomatic number is equal to 5. All the vertexes of the $H_{k3}$ graph (excluding both added $x$ and $y$ vertexes) correspond to the ordered pairs of the $H_{k2}$ graph's vertexes. Every ordered triple of the $H_{q1}$ graph's vertexes correspond to some $H_{k3}$ graph's vertex (fig. 9).

Let us execute the converting operation once more. We'll arrange table 5 (fig. 10) by the $F_{q3}$ matrix similarly like we arranged tables both 1 and 3.

The first two columns of table 6 are arranged by table 5 (similarly to tables both 2 and 4).

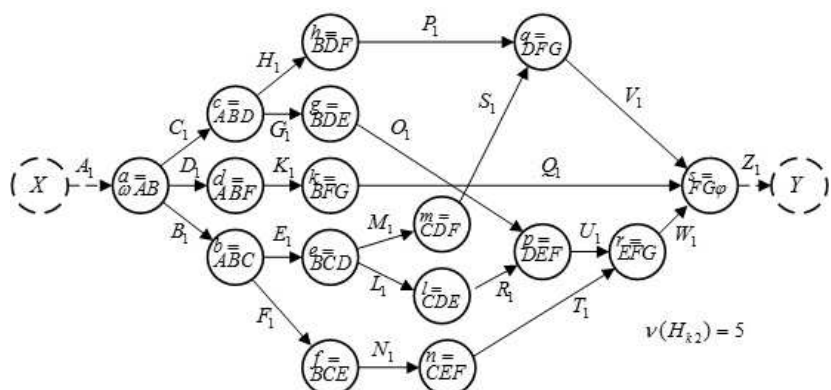

Fig. 9.

Table 5

| $N_{нач}$ | $q_i$ | $N_{кон}$ |
|---|---|---|
| x | $A_1$ | a |
| a | $B_1$ | b |
| a | $C_1$ | c |
| a | $D_1$ | d |
| b | $E_1$ | e |
| b | $F_1$ | f |
| c | $G_1$ | g |
| c | $H_1$ | h |
| d | $K_1$ | k |
| e | $L_1$ | l |
| e | $M_1$ | m |
| f | $N_1$ | n |
| g | $O_1$ | p |
| h | $P_1$ | q |
| k | $Q_1$ | r |
| l | $R_1$ | p |
| m | $S_1$ | q |
| n | $T_1$ | r |
| p | $U_1$ | r |
| q | $V_1$ | s |
| r | $W_1$ | s |
| s | $Z_1$ | y |

Table 6

| Vertexes of the $G_{φ3}$ graph | Pairs of vertexes of the $H_{φ3}$ graph | Pairs of the tripples of the vertexes of the $H_{φ2}$ graph | | Fours of the vertexes of the $H_{ψ1}$ graph |
|---|---|---|---|---|
| $A_1$ | xa | x | ωAB | xωAB |
| $B_1$ | ab | ωAB | ABC | ωABC |
| $C_1$ | ac | ωAB | ABD | ωABD |
| $D_1$ | ad | ωAB | ABF | ωABF |
| $E_1$ | be | ABC | BCD | ABCD |
| $F_1$ | bf | ABC | BCE | ABCE |
| $G_1$ | cg | ABD | BDE | ABDE |
| $H_1$ | ch | ABD | BDF | ABDF |
| $K_1$ | dk | ABF | BFG | ABFG |
| $L_1$ | el | BCD | CDE | BCDE |
| $M_1$ | em | BCD | CDF | BCDF |
| $N_1$ | fn | BCE | CEF | BCEF |
| $O_1$ | gp | BDE | DEF | BDEF |
| $P_1$ | hq | BDF | DFG | BDFG |
| $Q_1$ | ks | BFG | FGφ | BFGφ |
| $R_1$ | lp | CDE | DEF | CDEF |
| $S_1$ | mq | CDF | DFG | CDFG |
| $T_1$ | nr | CEF | EFG | CEFG |
| $U_1$ | pr | DEF | EFG | DEFG |
| $V_1$ | qs | DFG | FGφ | DFGφ |
| $W_1$ | rs | EFG | FGφ | EFGφ |
| $Z_1$ | sy | FGφ | y | FGφy |

Fig. 10

Using these tables, let's decode the ordered pairs of the $F_{q3}$ matrix's rows, which are written in the second column of table 6. For example, the $ab$ pair (the second line, the second column of table 6) includes both $a = \omega AB$ and $b = ABC$ elements. Therefore, in the second line of the third column of table 6 by this time the pairs of triples both $\omega AB$ and $ABC$ elements are written. All the third column of table 6 (fig. 11) is filled up in such a way.

Let's note that in all the lines of the third column of table 6 two last elements of the first triples are congruent to the first two elements of the second triples.

Therefore, every ordered pair of the triples corresponds to the ordered four of the $H_{q1}$ graph's vertexes, excluding both the first and the last ordered fours, each including two vertexes of the initial graph and two additional. One vertex was added to the $H_{q2}$ graph, and the second was added to the $H_{q3}$ graph. For example, the pair of both $\omega AB$ and $ABC$ triples correspond to the $\omega ABC$ four.

The ordered fours are written in the fourth column of table 6. The $R_{k3}$ matrix is arranged according to table 6 (fig. 11). The matrix is transformed into the $F_{k4}$ matrix (fig. 12) similarly to both $F_{k2}$ and $F_{k3}$ matrixes were arranged.

$F_{k3} =$

|   | x | a | b | c | d | e | f | g | h | k | l | m | n | p | q | r | s | y |
|---|---|---|---|---|---|---|---|---|---|---|---|---|---|---|---|---|---|---|
| x | \ | $A_1$ |   |   |   |   |   |   |   |   |   |   |   |   |   |   |   |   |
| a |   | \ | $B_1$ | $C_1$ | $D_1$ |   |   |   |   |   |   |   |   |   |   |   |   |   |
| b |   |   | \ |   |   | $E_1$ | $F_1$ |   |   |   |   |   |   |   |   |   |   |   |
| c |   |   |   | \ |   |   |   | $G_1$ | $H_1$ |   |   |   |   |   |   |   |   |   |
| d |   |   |   |   | \ |   |   |   |   | $K_1$ |   |   |   |   |   |   |   |   |
| e |   |   |   |   |   | \ |   |   |   |   | $L_1$ | $M_1$ |   |   |   |   |   |   |
| f |   |   |   |   |   |   | \ |   |   |   |   |   | $N_1$ |   |   |   |   |   |
| g |   |   |   |   |   |   |   | \ |   |   |   |   |   | $O_1$ |   |   |   |   |
| h |   |   |   |   |   |   |   |   | \ |   |   |   |   |   | $P_1$ |   |   |   |
| k |   |   |   |   |   |   |   |   |   | \ |   |   |   |   |   | $Q_1$ |   |   |
| l |   |   |   |   |   |   |   |   |   |   | \ |   |   | $R_1$ |   |   |   |   |
| m |   |   |   |   |   |   |   |   |   |   |   | \ |   |   | $S_1$ |   |   |   |
| n |   |   |   |   |   |   |   |   |   |   |   |   | \ |   |   | $T_1$ |   |   |
| p |   |   |   |   |   |   |   |   |   |   |   |   |   | \ |   | $U_1$ |   |   |
| q |   |   |   |   |   |   |   |   |   |   |   |   |   |   | \ |   | $V_1$ |   |
| r |   |   |   |   |   |   |   |   |   |   |   |   |   |   |   | \ | $W_1$ |   |
| s |   |   |   |   |   |   |   |   |   |   |   |   |   |   |   |   | \ | $Z_1$ |
| y |   |   |   |   |   |   |   |   |   |   |   |   |   |   |   |   |   | \ |

Fig. 11.

Then the $H_{k4}$ graph is constructed by this matrix. It appears to be the canonical graph. Its cyclomatic number is equal to 5. All the vertexes of the $H_{k4}$ graph (excluding the added both Γ and Ω vertexes) correspond to the ordered pairs of the $H_{k3}$ graph's vertexes. Every ordered pair of the $H_{k1}$ graph's vertexes corresponds to some $H_{k4}$ graph's vertex (fig. 13).

So, the straight converting allows constructing a consecutive row of the graphs, which vertexes correspond to the increasingly long orders of the initial graph's vertexes. In other words, the $H_j$ graph, received as a result of the consecutive straight converting of the $H_1$ graph, being done $(j-1)$ number of times (where $n = j$), possesses the following feature: any $H_j$ graph's vertex corresponds uniquely to every ordered $n$ of $H_1$ graph's vertexes.

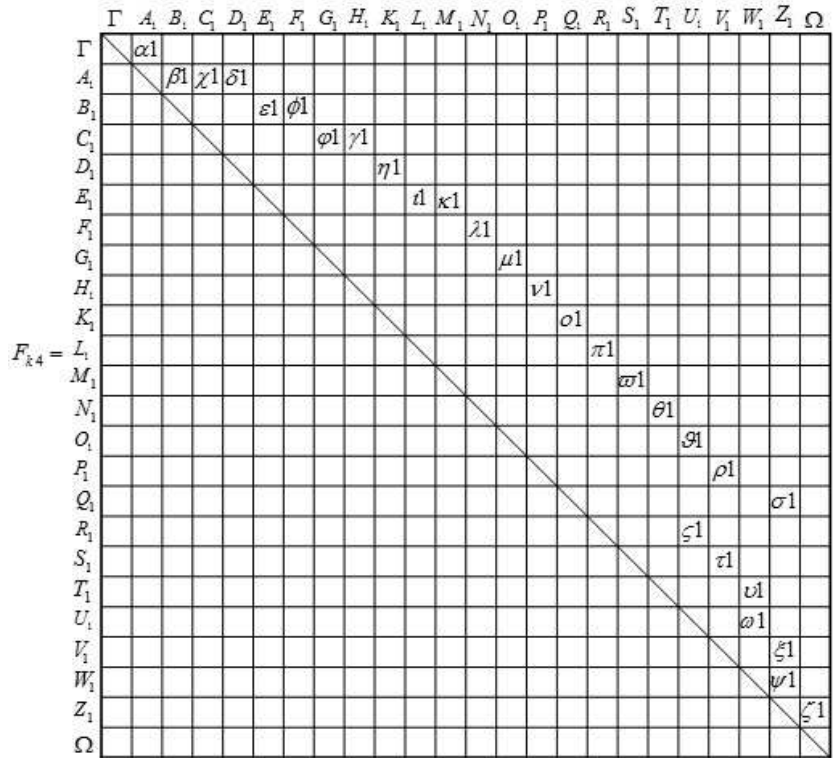

Fig. 12.

After the straight converting of the $H_1$ graph $(m-1)$ times ($m$ – is the order of the $F_1$ matrix) we'll get the $H_m$ graph, in which a part of the vertexes corresponds to all possibly ordered $m's$ of the vertexes of the $H_1$ graph.

An examination of such consecutive converting process allows doing the following evident conclusions:

1. The received $m's$ of the $H_1$ graph's vertexes corresponds to all the directed paths of the $(m-1)$ length, which are available in the $H_1$ graph.

2. Among the received paths there will be all the contours of not more than the $(m-1)$ length.

3. While accomplishing the converting $m$ times, we'll get all possible ordered $m's$, among which there will be contained all the $H_1$ graph's contours of not more than the $m$ length, including all the Hamilton circuits.

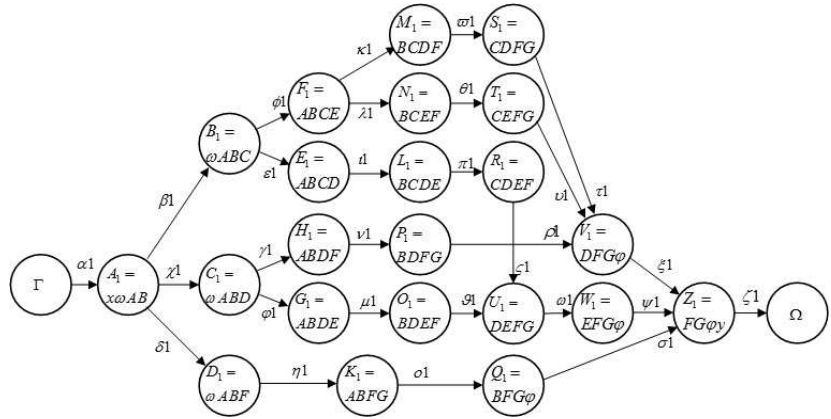

Fig. 13.

Thus, the straight converting allows constructing all possible both paths and contours of the given length in the initial graph by a rather simple approach. But this approach in cases, when the total number of both the paths and the contours in the $H_1$ graph is too large and the number of the $H_j$ graph's vertexes at the consecutive converting increases quickly, turns out to be inconvenient for searching either all or the part of the paths of the given form.

## 3 The interconnection between the cyclomatic number, the type of the converting and the increasing of the number of the vertexes

Let's show that the increasing of the number of the $H_j$ graph's vertexes, derivable by the consecutive converting is concerned with the cyclomatic number and the form of the straight converting (either the quasicanonical or the canonical).

Let the $H_j$ graph be received after the consecutive straight converting of the $H_j$ graph $(j-1)$ times.

The cyclomatic number of the $H_j$ graph will be:

$v(H_j) = m_j - n_j + 1$

Essentially for the $H_{(j-1)}$ graph:

$$v(H_{(j-1)}) = m_{(j-1)} - n_{(j-1)} + 1$$

From the content of the converting operation at the condition that the $H_1$ graph has both initial and final edges it follows that:

$$n_{(j+1)} = m_j + 2$$

And

$$n_j = m_{(j-1)} + 2$$

From which

$$m_j = n_{(j+1)} - 2$$

And

$$m_{(j-1)} = n_j - 2$$

From the expressions above we'll have:

$$n_{(j+1)} = n_j + v(H_j) + 1$$
$$n_j = n_{(j-1)} + v(H_{(j-1)}) + 1$$

Subtracting the latest expression from the last but one, we'll get:

$$\Delta n_{(j+1)} = \Delta n_j + \Delta v(H_j)$$

From the expression above it follows that if the given $H_1$ graph in the process of the consecutive straight converting generates only the canonical $H_{kj}$ graphs, for which, due to theorem 4 [2]: $\Delta v(H_j) = 0$.

Then the process of the consecutive straight converting of such $H_1$ graph is characterized by the following condition: $\Delta n_j = Const(j)$.

If the $H_1$ graph has both the initial and final edges, than:

$$\Delta n_j = v(H_1) - 1$$

In other words, the increasing of the number of the $H_{kj}$ graph's vertexes in comparison with the $H_{k(j-1)}$ graph does not depend on the $(j-1)$ number of the converting step, but depends only on the structure of the initial $H_1$ graph.

So we proved the following theorem:

### Theorem 9

If the $H_1$ graph in the process of its consecutive straight converting generates only the canonical $H_{kj}$ graphs at $(j = 1, 2, 3 ... M)$, then the numbers of the vertexes of these graphs, received step by step, are determined by the linear dependence from the converting operation's $(j-1)$ number, that is:

$$n_j = n_1 + \Delta n * (j - 1)$$

Relying on the conducted above argumentation also appears evident

## *Theorem 10*

If the $H_1$ graph in the process of its consecutive straight converting generates both the canonical and the quasicanonical graphs or only the quasicanonical graphs, then the numbers of the vertexes of these graphs, received step by step, are determined by the following expression:

$n_j = n_1 + \sum_{\xi=1}^{\xi=(j-1)} \Delta n_\xi$

Where: $\Delta n_\xi = \Delta n_{(\xi-1)} + \Delta v(H_{(\xi-1)})$

Where: $\xi = 1, 2, \dots (j-1), j = 1, 2, 3, \dots M$

For one's turn the $\Delta v(H_{(\xi-1)})$ values are determined by the structure of the $H_1$ graph.

From the formulations of theorems 9 and 10 it follows that all the directed graphs may be divided into two classes:

- The graphs, for which the cyclomatic number always appears to be the invariant of the straight converting on the arbitrary given above converting step.
- All the other directed graphs or graphs, for which the cyclomatic number on either the part or on all the steps of the converting, do not appear to be the invariant of the converting.

Let's proceed to reveal both the necessary and the sufficient characteristics of these classes of the directed graphs, depending on the concept of the path in the directed graph [3].

By the path in the $H(V, Q)$ graph we mean such a consecution of the $(q_1, q_2 \dots)$ edges that the exit of each previous edge coincides with the beginning of the following edge. The path is called simple if none edges meet twice and the compound otherwise. The path in which no vertex meets twice is called elementary. A contour is the finite path in which the final $v_k$ vertex coincides with the initial $v_1$ vertex. A contour is called elementary if all its vertexes are different, except coincident both the initial and the final vertexes.

Let bring in one more feature, which allows setting apart one path from the other in the canonical graphs. Let us determine this feature with the help of the degrees of the $\rho(v_h)$ vertexes of the $H$ graph, through which this path goes along. Let's denote the arbitrary path as $\xi$ and arrange the function for this path:

$\rho_{i_\xi}(v_h) = \sum_1^{i_\xi} \rho_\chi(v_h)$

Where: $\chi = 1, 2, 3 \dots i_\xi$ and $i_\xi = 1, 2, 3, \dots l_\xi \dots n_\xi$ – the numbers of the $v_h \in \xi$ vertexes.

Since we examined only the graphs, which correspond to the finite real processes, then every contour (due to the accepted admission) is included into either that or other path only the finite number of times. Therefore, every vertex may meet with the contour only the finite number of times and every time it will be examined as a new vertex.

The $\rho_{i_\xi}(v_h)$ functions may be divided into two $A_p$ and $B_p$ classes, which we'll define the following way.

The $\rho_{i_\xi}^{(A)}(v_h) \in A_p$ function on the $1 \leq i_\xi \leq i_\xi^*$ interval represents the never-increasing function of $i_\xi$, and on the $i_\xi^* \leq i_\xi \leq n_\xi$ interval – the never-decreasing function of $i_\xi$. Let suppose that all the functions, unsatisfying this condition, are: $\rho_{i_\xi}^{(B)}(v_h) \in B_p$.

From the definition of both $A_p$ and $B_p$ classes, and the $\rho_{i_\xi}(v_h)$ functions it is evident that functions, definitely subordinated to the unified rule for all the functions, belong to the $A_p$ class. Let us agree to denote such functions as $\rho_{i_\xi}^{(A)}(v_h) \in A_p$ and the corresponding paths as the holonomic paths.

The $\rho_{i_\xi}(v_h)$ functions, subordinated to any arbitrary rule at the $i_\xi$ functional changes, excluding the rule, determined for the $A_p$ class, belong to the $B_p$ class. Let us agree to denote such functions as $\rho_{i_\xi}^{(B)}(v_h) \in B_p$ and corresponding paths as the heteronomous paths.

Now let us prove the theorem on the converting invariant.

### Theorem 11

A cyclomatic number appears to be the invariant of the straight converting of the $H_j$ graphs at any $j = 1, 2, 3, \ldots M$ value if and only if the $H_1$ graph is canonical and all the paths in the $H_1$ graph are the holonomic paths.

In other words, the consecutive straight converting of the initial canonical $H_1$ graph with the holonomic paths generates the $H_j$ graphs with the same cyclomatic number as the initial graph has.

### Proof

From theorem 4 [2] it follows that the cyclomatic $\nu(H_j)$ number of the $H_j$ graph at the straight converting of the $H_j$ graph into the $H_{(j+1)}$ graph does not change if and only if the $H_j$ graph is a canonical one, in other words, if all the $H_j$ graph's vertexes are simple ones (the simple vertexes also include the elementary vertexes). If on some $(j-1)$ step of the converting in the $H_j$ graph the complicated vertexes appear, then on the next $j$ step of converting a complicated vertex generates at least one independent cycle, which increases the cyclomatic number.

Therefore, for the demonstration of the theorem, it is enough to show that the complicated vertexes arise in the $H_j$ graph in no $(j-1)$ steps of the $H_1$ graph's converting if and only if all the paths in the $H_1$ graph are the holonomic ones.

Let us examine the arbitrary $\xi$ path, which may include the contours and, therefore, may consist both recurring and non-recurring both edges and vertexes. For the proof of the theorem it is enough to consider any contour as the iterative at least once.

In every simple $v_h$ vertex of the $H_1$ graph, excluding both initial and final vertexes, we have either the divergent fan of the edges at $\rho(v_h) < 0$ (fig. 14a), or the converging fan of the edges at $\rho(v_h) > 0$ (fig. 14b). The vertexes for which $\rho(v_h) > 0$ (fig. 14b) we'll denote as $v_h^{(+)}$ (a positive vertex) and for which $\rho(v_h) < 0$ (fig. 14a) we'll denote as $v_h^{(-)}$ (a negative vertex). The elementary vertexes we'll denote as $v_h^{(0)}$.

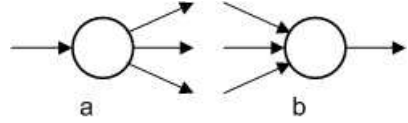

Fig. 14.

Let's also bring in the notations for the $l_{xy}$ intervals along the $\xi_j$ path between the vertexes of the different types. The notations of the corresponding intervals are listed in the cells of table 7 (fig. 15) at the left of the bottom. The other necessary explanations to table 7 will be done afterwards.

Table 7

| Intervals between the vertexes of different types and their characteristics | | | | | | |
|---|---|---|---|---|---|---|
| Initial vertex of the interval \ Final vertex of the interval | $v_{h(i_e,j)}^{(+)}$ | | $v_{h(i_e,j)}^{(0)}$ | | $v_{h(i_e,j)}^{(-)}$ | |
| $v_{h(i_e,j)}^{(+)}$ | $l_{11}$ | $A_p; B_p$ | $l_{12}$ | $A_p; B_p$ | $l_{13}$ | $A_p; B_p$ |
|  | 0 | (−) | − | − | 1 | (−) |
| $v_{h(i_e,j)}^{(0)}$ | $l_{21}$ | $A_p; B_p$ | $l_{22}$ | $A_p; B_p$ | $l_{23}$ | $A_p; B_p$ |
|  | − | − | − | − | − | − |
| $v_{h(i_e,j)}^{(-)}$ | $l_{31}$ | $B_p$ | $l_{32}$ | $A_p; B_p$ | $l_{33}$ | $A_p; B_p$ |
|  | −1 | (+) | − | − | 0 | (−) |

Fig. 15.

The following statements are evident:

- If on the $i_\xi^{(beg)} \leq i_\xi \leq i_\xi^*$ interval of the $\xi$ path every next vertex along this path is either positive or elementary, then on this interval the $\rho_{i_\xi}(v_h)$ function appears to be the never-increasing one;

- If on the $i_\xi^* \leq i_\xi \leq i_\xi^{(fin)}$ interval of the $\xi$ path every next vertex along this interval is either negative or elementary, then on this interval the $\rho_{i_\xi}(v_h)$ function appears the never-decreasing one.

It follows that

### Lemma 1

The $\xi_1 \in A_p$ path is a holonomic one if and only if any positive vertex along the $\xi_1$ path has the lesser ordinal number then any of the negative vertexes.

Let us prove two statements.

1) If the $\xi_1 \in A_p$ path, then at no $(j-1)$ step of the straight converting on the $\xi_j$ path, which is the result of this step, the complicated vertex could be generated.

2) In order to complicated vertex be generated on the $\xi_j$ path, which is the result of the $(j-1)$ step of converting of the $H_1$ graph, containing the corresponding $\xi_1$ path, it is necessary and sufficient for the $\xi_1$ path to be the heteronomous one.

Let the generated after the $(j-1)$ converting step the $v_{h(i_\xi,j)}^{(+)}$ vertex has along the $\xi_j$ path the ordinal $i_\xi$ number. The distance from this vertex to the vertex, which is accepted as the initial vertex of the $\xi_j$ path, makes up the $(i_\xi - 1)$ edges (fig. 16a).

At converting the $H_j$ graph into the $H_{(j+1)}$ graph instead of the last edge of this interval, having the ordinal $(i_\xi - 1)$ number, a new $v_{h(i_\xi-1),(j+1)}^{(+)}$ vertex is generated, which has the same number (fig. 16b). In fig. 16a a path is shown, which is received after the $(j-1)$ step of converting of some $\xi_1$ path, which enters the initial $H_1$ graph. The $\xi_{(j+1)}$ path, received as a result of converting of the $\xi_j$ path, is presented in fig. 16b. In fig. 16a the edges, which were transformed into the vertexes at graph's converting (fig. 16b), are marked with the dots in fig. 16a.

So, at every converting step of the $H_j$ graph instead of every positive $H_j$ graph's vertex in the $H_{(j+1)}$ graph a corresponding also positive vertex is generated, but having the ordinal number less per unit.

In case if the $H_1$ graph has the initial edge, then at every converting step the $H_j$ graph is increased, as it was shown before, with the new initial edge. It should be done even if the $H_1$ graph consists of a single unlocked path. Then the ordinal numbers of the positive $H_{(j+1)}$ graph's vertexes corresponding to the positive $H_j$ graph's vertexes do not change.

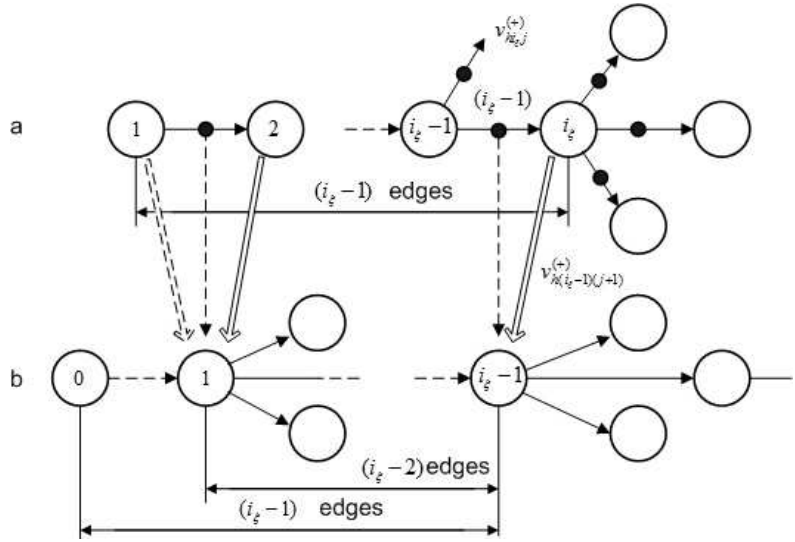

Fig. 16.

Then, let some negative $v^{(-)}_{h(i_\xi j)}$ vertex, generated after the $(j-1)$ converting step of the $H_1$ graph has along the $\xi_j$ path the ordinal $i_\xi$ number. The distance to this vertex from the vertex, which is accepted as the initial vertex of the $\xi_j$ path, has, as it was in the previous case, the $(i_\xi - 1)$ edges (fig. 17a). While converting the $H_j$ graph into the $H_{(j+1)}$ graph instead of the edge, having the ordinal $i_\xi$ number, the negative $v^{(-)}_{h(i_\xi(j+1))}$ vertex, having the ordinal $i_\xi$ number, will be generated (fig. 17b). Even if the $H_1$ graph had the initial edge, then the ordinal number of the newly generated $v^{(-)}_{h(i_\xi+1)(j+1)}$ vertex will be per unit more, then the number of corresponding to it the $v^{(-)}_{h(i_\xi,j)}$ vertex, that is $(i_\xi + 1)$.

The $\xi_j$ path is shown in fig. 17a. It is received after the $(j-1)$ converting step of some $\xi_j$ path, entering into the initial $H_1$ graph. The $\xi_{(j+1)}$ path, received as a result of the $\xi_j$ path converting, is shown in fig. 17b. The edges corresponding to the vertexes in fig. 17b are presented with the dots in fig. 17a.

So, we have determined the dependence of the ordinal numbers of the vertexes of the different types from the number of the step of the $H_j$ graph's converting. It allows estimating the modification of the $l_{xy}$ intervals at the graph's converting (tab. 7, fig. 15). Let us estimate this modification at one step of converting and fill in table 7 in the following order:

On the left at the top – the symbol of the $l_{xy}$ interval;

On the right at the top – the classes of the $\xi$ paths, to which the $l_{xy}$ interval may belong;

On the left at the bottom – the increase of the length of the corresponding $l_{xy}$ interval at the converting of the $H_j$ graph into the $H_{(j+1)}$ graph;

On the right at the bottom – with sign (+) or (−) accordingly is marked either the possibility or the impossibility of the generation of the complicated vertex on changing of the interval's length to zero at the graph converting.

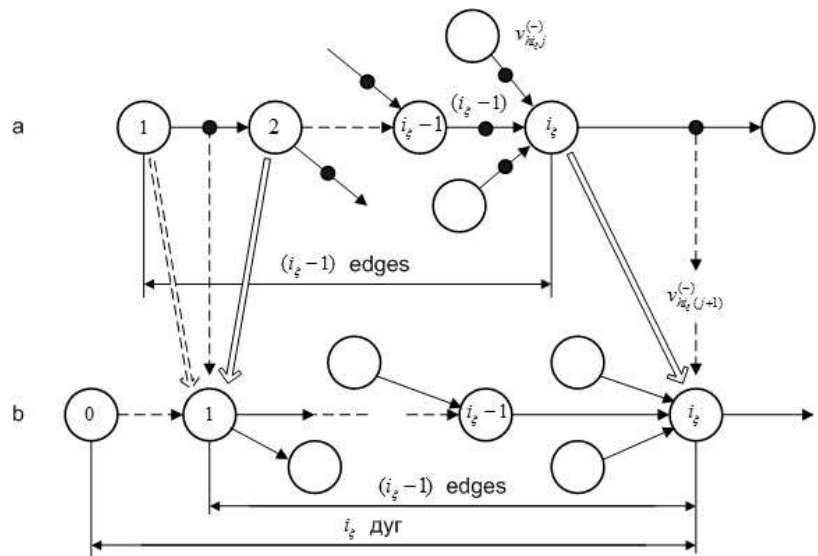

Fig. 17

The following conditions of filling in table 7 (fig 15) are evident from the argumentations above:

1. The $l_{11}$ interval may belong to the $\xi_j$ paths of both $A_p$ and $B_p$ classes. At the $H_j$ graph's converting the $l_{11}$ interval generates in the $H_{(j+1)}$ graph the interval of the same type ($l_{11}$) and length. So the modification of the ordinal numbers of both the initial and the final vertexes of both $l_{11j}$ and $l_{11(j+1)}$ intervals are equal, in other words, $\Delta l_{11(j+1)} = 0$. Therefore a complicated vertex instead of the $l_{11j}$ interval could be generated at no converting step.

2. The $l_{12}$, $l_{21}$, $l_{23}$ and $l_{32}$ intervals have the elementary vertexes on one of the endings. The $l_{22}$ interval has the elementary vertexes on the both endings. As a result, the complicated vertexes could not be generated instead of these intervals at the graph's converting.

3. The $l_{13}$ interval may belong to the $\xi_j$ paths of both $A_p$ and $B_p$ classes. This interval at the graph's converting generates the interval of the same type, but with the length more per unit. Therefore, no complicated vertexes could be generated instead of the $l_{13}$ interval at any converting step.

4. An $l_{33}$ interval may belong to the $\xi_j$ paths of both $A_p$ and $B_p$ classes. The $l_{33}$ interval generates in the $H_{(j+1)}$ graph the interval of the same type and length at every converting step. Therefore, the complicated vertex could be generated on no steps of converting of the $H_j$ graph into the $H_{(j+1)}$ graph.

5. An $l_{31}$ interval may belong only to the $\xi_j$ path of the $B_p$ class. It generates the $l_{31(j+1)}$ interval of the same type, but having the length less per unit than the $l_{31}$ interval, at the converting of the $H_j$ graph. As a result the $l_{31(j+1)}$ interval is generated instead of the $l_{31}$ interval at some finite $j$ step of converting. It has a zero length. It means the appearance of the complicated vertex (fig. 18).

A process of the transformation of the $l_{31}$ interval into the complicated vertex at the $H_j$ graph's converting is presented in fig. 18 (fig. 18a, 18b, 18c – the converting steps $((j-2), (j-1), ...)$. This vertex on the $(j+1)$ converting step generates the four independent cycles (fig. 18d).

Both formulated above statements (1) and (2) evidently follow from table 7.

So, for the cyclomatic number of the $H_{k1}$ graph to be the invariant of converting (at any number of the steps), it is necessary and sufficient for the $H_{k1}$ graph to contain none of the heteronomous paths.

The last statement proves the theorem.

Theorem 11 has the evident

### *Corollary*

A highest possible $j_{max}$ of the converting steps at $v(H_{kj}) = Const(j)$, which the $H_1$ graph allows, is equal to the length of the shortest of all, belonging to the $H_1$ graph, the intervals of the $l_{31}$ type, in other words, it is equal to $\min_{\{\xi_1\}} l_{31}$, where all $\xi_1 \in H_1$.

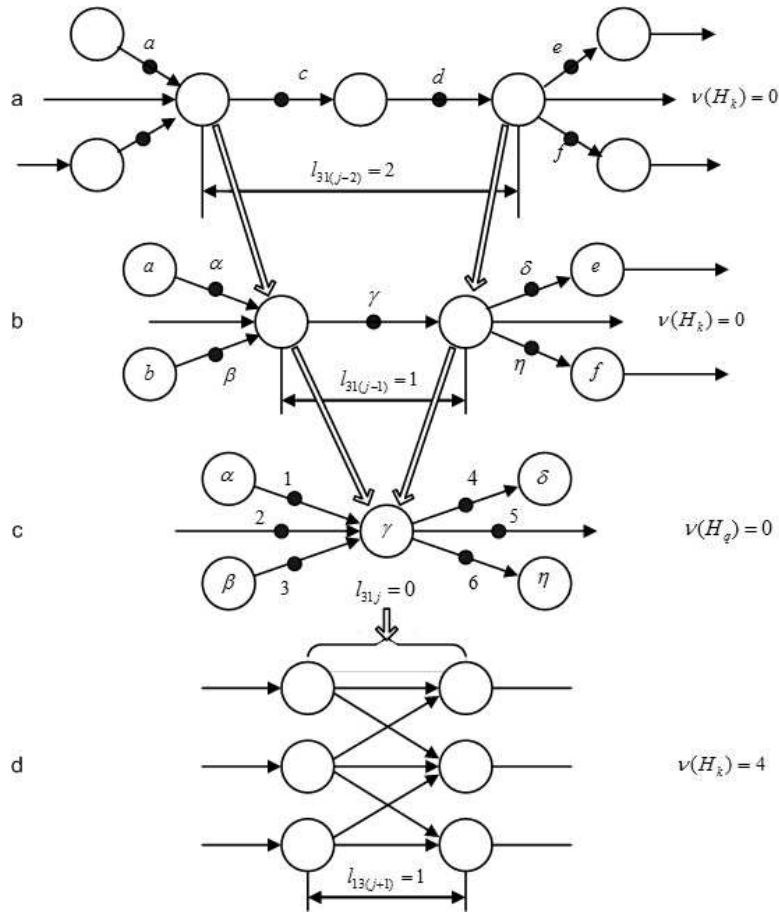

Fig. 18

## 4 The cyclomatic number and graph's classes

The $v(H)$ cyclomatic number appears to be one of the features of the $H$ graph's complexity. Therefore, from a certain point of view, we can say that if the cyclomatic number appears to be the converting invariant, then the complexity of the $H_j$ graphs, received as a result of converting, does not differ from the complexity of the initial $H_1$ graph.

So, we find out that some graphs permit the invariant (in the terms of the cyclomatic number) converting; other graphs do not permit such a converting. Let us show that all the concerned directed graphs may be divided into three classes:

$H^{(1)}$ – the graphs, which permit the consecutive invariant converting any number of times;

$H^{(2)}$ – the graphs, for which the cyclomatic number, beginning from some $j_{min} \geq 0$ step of the converting, is not the invariant of the converting.

But for such graphs we can indicate another, may be very big, but the finite $j_{kp}$ number of the converting steps, after the accomplishment of which, the cyclomatic number also begins and then continues to be the invariant of the converting, but it does not matter how many times the next converting is accomplished;

$H^{(3)}$ – the graphs, for which, after already accomplished the arbitrarily large number of $j = M$ steps of the preliminary converting, we can indicate the $j > M$ steps of the converting, at the accomplishment of which the increasing of the $H_{(j+1)}$ graph's cyclomatic number will again occur.

For the definition of both the necessary and the sufficient conditions of any of three classes let's examine more thoroughly the concept of the contour in the directed graph. The definition of the graph's contour, examined separately from the graph, which it is entering, was given above.

Every contour is arranged from some $Q_1 \in Q$ subset, where $Q$ – is the set of all the $H$ graph's edges. A $Q_1$ subset must be connected with the $Q_2$ subset of the other $H$ graph's edges. The edges, connecting the contour with the $Q_2$ subset of the other graph's edges, may be, concerning the contour, either the "entry" (fig. 19) or the "exit" (fig. 20).

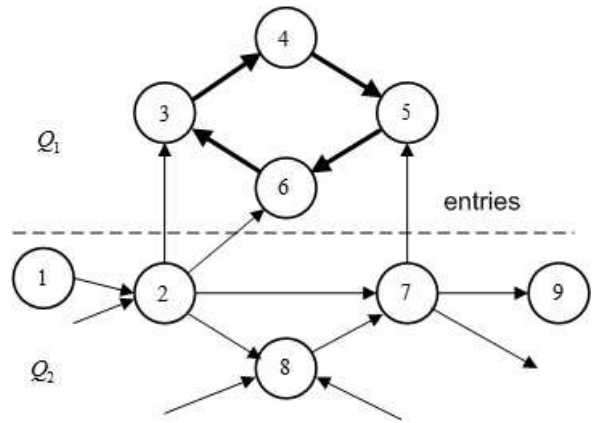

Fig. 19.

It is evident, that the real processes correspond only to such contours, which have at least one "entry" and at least one "exit" (fig. 21). The contours are marked by the thick lines in fig. 19, 20 and 21.

If the contour has not a single "exit", then such contour practically is isolated from the process. Let's call such contour as the contour-trap or the deadlock contour.

If the contour has not a single "entry", then such contour also is isolated from the remaining process. Let's agree to call such contour as the fictitious.

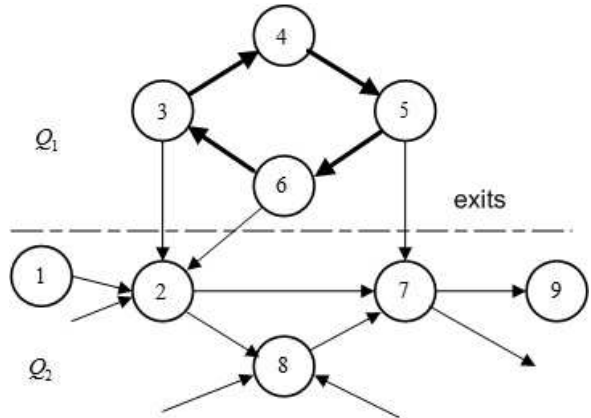

Fig. 20.

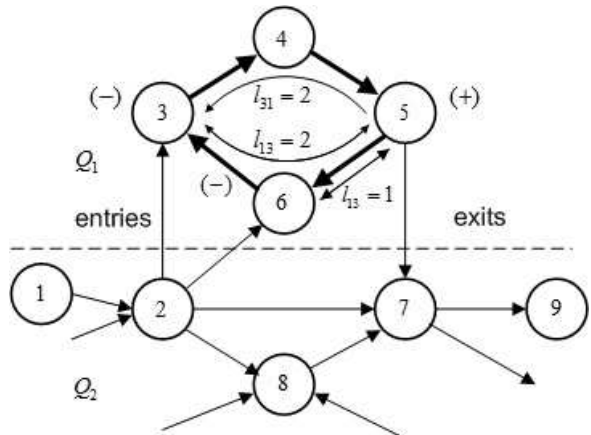

Fig. 21.

It is clear that neither the deadlock contour nor the fictitious contour could be examined as a part of the real process' model because with their help we could not interpret the seriated part of the elements of the process.

In this connection we'll examine only such contours, which have both the "entries" and the "exits". Let us agree that in all cases by the word contour we'll denote the contours, which have both the "entries" and the "exits". The exclusion is the graph, consisting of the elementary contour.

Let's prove three theorems, which define the necessary and sufficient features of the listed above three graph's classes.

## Theorem 12

A $H_{k1}$ graph belongs to the $H^{(1)}$ class if and only if it does not contain either the contours or the intervals of the $l_{31}$ type.

*Proof*

Accordingly to theorem 11 the $H_{k1}$ graph permits the invariant converting if and only if all the paths in the graph are the holonomic ones. A path is not holonomic (is heteronomous) if and only if it has at least one interval of the $l_{31}$ type.

A presence of the contours is equally matched to the presence of the intervals of the $l_{31}$ type in the graph. Indeed, if the contour does not contain the complicated vertex, then, based on the above argumentations (fig. 19, 20 and 21), the contour, having at least one both "entry" and "exit", should have at least one interval of the $l_{31}$ type and one interval of the $l_{13}$ type (fig. 21).

Thus, the $H_{k1}$ graph belongs to the $H^{(1)}$ class if and only if it does not contain the contours and the intervals of the $l_{31}$ type. The theorem is proved.

Let's denote the graphs of the $H^{(1)}$ class as a class of the holonomic graphs.

Now theorem 12 may be formulated as follows: a cyclomatic number is the regular invariant of the consecutive straight converting of the $H_{k1}$ graph independently from the number of the converting steps if and only if the $H_{k1}$ graph is the holonomic one, in other words, if the $H_{k1}$ graph contains only the holonomic paths.

*Remark*

It may be shown that the graphs with the isolated (both the fictitious contours and the contour-traps) contours, but without the heteronomous paths also belong to the holonomic class of the graphs. But in this work such graphs were not examined.

*Theorem 13*

A $H_1$ graph belongs to the $H^{(2)}$ class if and only if it contains the intervals of the $l_{31}$ type, but does not contain the contours.

*Proof*

If the $H_1$ graph does not contain the contours, then each of the intervals of the $l_{31}$ type in existence has the finite length. As far as the graph is being converted, the intervals of the $l_{31}$ type also generate the intervals of the $l_{31}$ type, but of a less length. As a result of the process in the graph each of the smallest intervals of the $l_{31}$ type turns to zero and at that generates the complicated vertex. The complicated vertexes at the next converting step generate each some finite number of the cycles, increasing at that the cyclomatic number of either $H_{kj}$ or $H_{qj}$ graph.

Each of new cycles consists of the intervals of the $l_{13}$ type (fig. 18), which, at further converting, do not influence the graph's cyclomatic number.

Since the number of the paths in the graph, the number of the intervals of the $l_{31}$ type and these lengths are finite, so some $j_{kp}$ number of the converting steps exists, after the fulfillment of which in the $H_{(j_{kp}+1)}$ graph there would remain not a single interval of the $l_{31}$ type and not a single complicated vertex.

From this converting step the cyclomatic $\nu(H_{kj})$ number of the $H_{kj}$ graph, where $j \geq j_{kp} + 1$, becomes the converting invariant.

So, as the $H_1$ graph belonged to the $H^{(2)}$ class, it is sufficient for it to contain the intervals of the $l_{31}$ type, but not to contain the contours. Let's show that this condition also serves as the necessary condition.

Let the graph have at least one contour (fig. 21).

This contour is examined separately in fig. 22a. If the contour has at least one "entry" and at least one "exit", then it has at least one interval of the $l_{31}$ type. The $l_{31} = 2$ interval connects both the $v_3$, $v_4$ and $v_5$ vertexes (fig. 22a).

The contour (fig. 22a) is locked by the $l_{13} = 2$ interval. A contour "$bcdeb$" (fig. 22b) is received at the converting from contour "34653" (fig. 22a). A converting leads to the decreasing of the interval of the $l_{31}$ type (fig. 22b), and then it turns to zero (fig. 22c). At that the complicated vertex is generated (fig. 22c – $v_\gamma$ vertex).

At the next converting step (fig. 23) the complicated vertex disintegrates into the series of the simple vertexes. At that the contour remains (fig. 23 – "34653") and the cycle appears (fig. 22 – "12341"). In the new contour the $l_{31}$ interval has the "3" length and the $l_{13}$ interval – the "1" length.

So, at the converting of any contour, the entering interval of the $l_{31}$ type turns to zero and generates the complicated vertex. Then the complicated vertex generates at least one cycle and one contour. The contour, in its turn, by all means contains the interval of the $l_{31}$ type.

Thus, the consecutive converting of the $H_j$ graph, which contains at least one contour, always leads to the appearance of at least one complicated vertex or the interval of the $l_{31}$ type. So it could not generate the graph of $H^{(1)}$ type.

Well, the absence of the contours appears to be the necessary condition of the $H_1$ graph belonging to the $H^{(2)}$ class. But for the $H_1$ graph it is necessary to contain at least one interval of the $l_{31}$ type in order not to belong to the $H^{(1)}$ class.

Once and for all: in order to belong to the $H^{(2)}$ class it is necessary and sufficient for the $H_1$ graph to contain the intervals of the $l_{31}$ type, but not the contours.

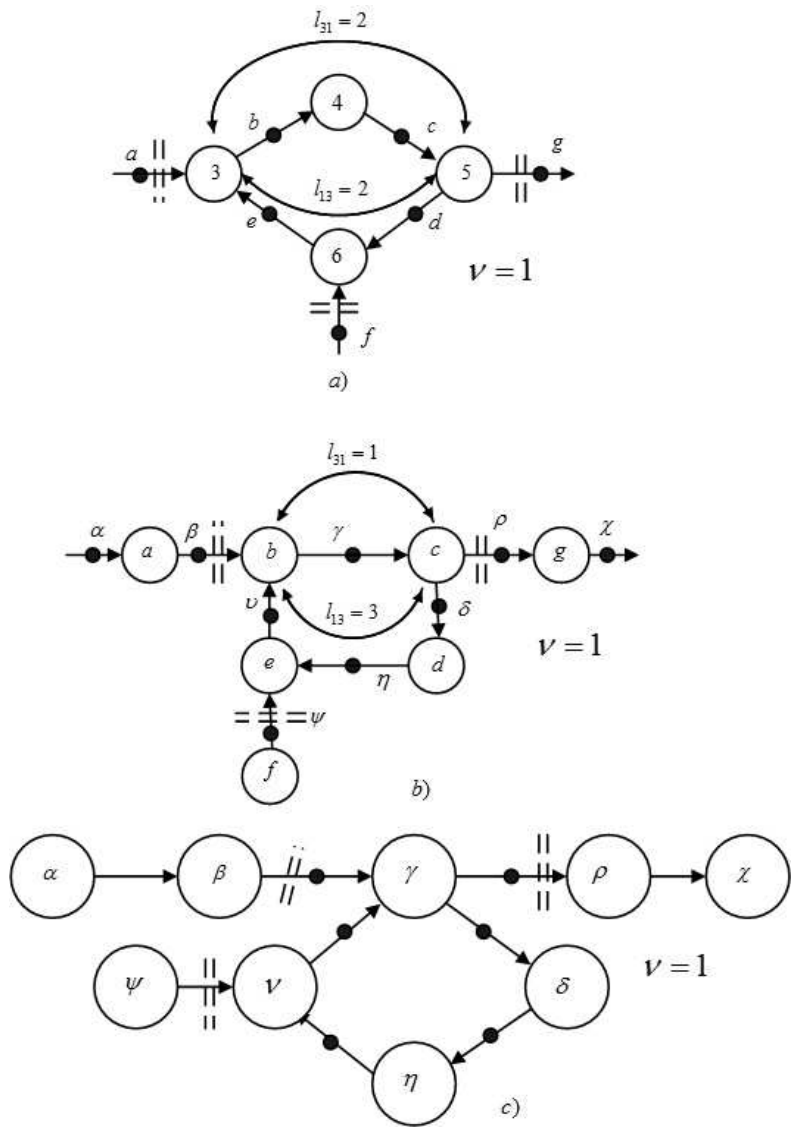

Fig. 22.

Let us denote the graphs of this class as the bounded-heteronomous graphs.

So, theorem 13 may be now formulated in the following way.

If the $H_1$ graph contains the heteronomous paths, but does not contain the contours, then it appears to be bounded-heteronomous and has, by the number of the steps of the consecutive straight converting, the heteronomous $(j_{kp})$ bound, up to which the cyclomatic number does not appear to be a regular invariant of the $H_j$ graph's converting. On reaching the heteronomous bound, the $H_{kj}$ graph becomes holonomic, and its cyclomatic number becomes a regular invariant of converting independently from the number of further steps of converting.

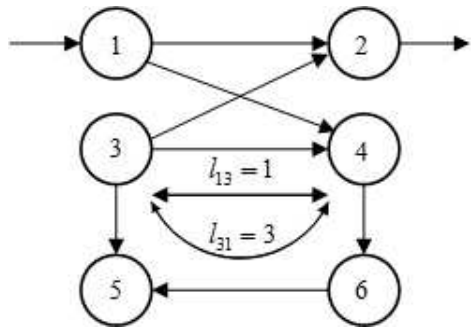

Fig. 23.

In accordance with the names, accepted for both $H^{(1)}$ and $H^{(2)}$ graph's classes, it is naturally to denote the graphs of the $H^{(3)}$ class as progressive-heteronomous graphs.

In connection with theorem 13 the following theorem 14 is evident.

## Theorem 14

So as the $H_1$ graph was progressive-heteronomous, it is necessary and sufficient for it to contain at least one contour.

A progressive-heteronomous graph, according to the $j$ number of the steps ($j$ may be arbitrarily large) has not the heteronomous bound on reaching of which the cyclomatic number becomes a regular invariant of converting.

In conclusion let's do a brief remark about the possible method of the establishment of the graph's belonging to either that or other class.

A graph's belonging to the progressive-heteronomous class is determined by the presence of at least one non-isolated contour.

A graph's belonging to the bounded-heteronomous class may be determined easily. It is sufficient to apply to the $H_k$ graph Faulks algorithm [4], which allows fixing the transit paths of the definite length between the given vertexes in the graph.

At that it is sufficient with the help of Faulks algorithm [4] to fix in the graph at least one path with the length not more than a unit, which is going from the negative vertex to the positive vertex, in other words, the heteronomous path.

If the $H_k$ graph is neither progressive-heteronomous nor bounded-heteronomous, then it is holonomic.

An examining of the operation of the straight converting allows formulating the following conclusions on both the possibility and the conditions of the multiple consecutive transformations of the directed edge graphs into the directed vertex graphs.

1. If the adjacency vertex matrix of the arbitrary directed $H_1$ graph meets the requirements, which are produced to the $F$ matrix – the adjacency matrix of the vertexes of the edge graph, then the $H_1$ graph can be subjected to the operation of the consecutive straight converting any number of times. At that all the graphs are divided as regards to the operation of the consecutive straight converting into three classes:

— A homonomic

— A bounded-heteronomous

— A progressive-heteronomous

2. For the homonomic graphs the cyclomatic number is a regular invariant of the converting independently from the number of the consecutive converting steps. Owing to this fact, the number of the $H_j$ graph's vertexes, received from the initial $H_1$ graph as a result of its consecutive converting, has the linear dependence from the number of the steps.

3. The bounded-heteronomous graphs have the heteronomous bound according to the number of the steps of the consecutive converting. Before the achievement of this bound, the cyclomatic number does not appear to be a regular invariant of converting, but on reaching this bound, as a result of the regular converting step, the homonomic graph is generated, and the cyclomatic number becomes the regular invariant of converting independently from the number of the steps of the following consecutive converting.

4. The progressive-heteronomous graphs do not have the heteronomous bound according to the number of the steps of the consecutive converting. As a result, the cyclomatic number of the progressive-heteronomous graph never becomes the regular invariant of converting even if the number of the steps will be enormously large.

## 5 Acknowledgments

Many thanks for help to improve this article to my friend Olga Volgina, who for several years had reformed my English. Special thanks for the members of my family who undergone all the difficulties side by side with me and who encouraged me in my work. Thanks to Michael Trofimov for his help in performing the material to arXiv.

Also I will be very grateful to those readers, who will find and send me a word about the uncovered misprints in order to improve the text. The matter is that on the base of the proved theorems were made the polynomial algorithms. It was also proved that the algorithms are local. And, finally, these algorithms were brought to the view of the programs.

## 6  Some designations

$F$ − An adjacency matrix of the $H$ graph's vertexes;
$R$ − An adjacency matrix of either the $G$ graph's vertexes or the $H$ graph's edges;
$L$ − An adjacency matrix of the initial $G$ graph's vertexes;
$J$ − An operation of the straight converting;
$D$ − An operation of the reverse converting;
$\nu(H)$ - A cyclomatic number;
$\xi_i$ − An arbitrary path;
$v_h^{(+)}$ − A positive vertex;
$v_h^{(-)}$ − A negative vertex;